\documentclass[fleqn]{mat01}
\usepackage{times,mathtimy,amssymb,latexsym}%,krepsf,rangecite}
\begin{document}

\setcounter{page}{267} \firstpage{267}

\def\dim{{\rm dim}}

\def\dim{\mathop {\rm dim}}
\def\emdim{\mathop {\rm emdim}}

\newtheorem{theore}{Theorem}
\renewcommand\thetheore{\arabic{section}.\arabic{theore}}
\newtheorem{theor}{\bf Theorem}
\newtheorem{lem}[theor]{\it Lemma}
\newtheorem{rem}[theor]{Remark}

\title{Maximally differential ideals in regular local rings}

\markboth{Alok Kumar Maloo}{Maximally differential ideals in
regular local rings}

\author{ALOK KUMAR MALOO}

\address{Department of Mathematics and Statistics, Indian
Institute of Technology, Kanpur~208~016, India\\
\noindent E-mail: akmaloo@iitk.ac.in}

\volume{116}

\mon{August}

\parts{3}

\pubyear{2006}

\Date{MS received 16 January 2006}

\begin{abstract}
It is shown that if $A$ is a regular local ring and $I$ is a
maximally differential ideal in $A$, then $I$ is generated by an
$A$-sequence.
\end{abstract}

\keyword{Regular rings; differential ideals; maximally
differential ideals.}

\maketitle

\section{Introduction}

Let $(A,\mathfrak{m})$ be a Noetherian local ring and let $I$ be
an ideal in $A$. Assume that $I$ is maximally differential under a
set of derivations of $A$. Then $A/I$ contains a field $k$. If $k$
is of characteristic $p>0$ then the structure of $A/I$ is given by
Harper's theorem \cite{lh} (see \cite{akm1} for a proof).

Assume now that $A/I$ contains a field of characteristic zero.
Then so does $A$ and by Lemma~1.3 of \cite{l}, $I$ is a prime
ideal. If $A$ is a $G$-ring (see p.~256 of \cite{hm} for
definition) then, by \cite{h}, $A/I$ is regular. However, in
general, $A/I$ is not regular (see \cite{l}, \cite{gl} or
\cite{akm3} for counter-examples). One notes that in all these
counter-examples, the ring itself was nonregular.

Therefore one may ask the following question: If $A$ is regular
then what are the properties of $A/I$? In particular, is $A/I$
regular? It is trivial to prove that $A/I$ is regular if
$\dim(A)\leq 1$ (note that $I$ is a prime ideal). If $\dim(A)\geq
2$, the problem appears to be hard and we have not been able to
solve it even for $\dim(A)=2$. However, what we show in this note
is that if $A$ is regular then $I$ is generated by an $A$-regular
sequence. Therefore $A/I$ is, in fact, a complete intersection.

We prove this result without any condition on the characteristic
of the ring. Moreover, $I$ may either be maximally differential
under a set of derivations or a set of higher derivations, i.e,
Hasse--Schmidt derivations.

We also extend our result (Theorem~4 of \cite{akm1}) about the
structure of a maximally differential ideal in positive
characteristic to unequal characteristic case.

\section{Results}

By a ring we mean a commutative ring with unity.

Let $A$ be a ring. We first recall a few definitions.

A derivation of $A$ is an additive endomorphism $d$ of $A$ such
that for all $a,b \in A$, $d(ab)=ad(b)+bd(a)$.

A higher derivation or a Hasse--Schmidt derivation $D$ of $A$ is a
sequence $(D_0,D_1,$ $D_2,\dots)$ of additive endomorphisms
$D_n$'s of $A$ such that $D_0$ is the identity of $A$ and for all
$n\geq 0$ and $a,b\in A$, $D_n(ab) = {\sum_{i+j=n}}D_i(a) D_j(b)$.

Let $d$  be a derivation of $A$. An ideal $I$ of $A$ is said to be
$d$-differential if $d(I)\subseteq I$. Similarly if
$D=(D_0,D_1,D_2,\dots)$ is a higher derivation of $A$ then an
ideal $I$ of $A$ is said to be $D$-differential if
$D_n(I)\subseteq I$ for all $n\geq 0$.

Let $\frak{D}$ be a set of derivations (or a set of higher
derivations) of $A$. An ideal $I$ of $A$ is said to be
$\mathfrak{D}$-differential if $I$ is $D$-differential for all
$D\in \mathfrak{D}$. An ideal $I$ of $A$ is called a maximally
$\mathfrak{D}$-differential ideal if it is a proper
$\mathfrak{D}$-differential ideal and for every ideal $J$ of $A$
with $I\subset J \subset A$, $J$ is not
$\mathfrak{D}$-differential.

If $A$ is local and $\mathfrak{D}$ is a set of derivations (or a
set of higher derivations) of $A$, then $A$ has a unique maximally
$\mathfrak{D}$-differential ideal $I$, and $I$ contains all proper
$\mathfrak{D}$-differential ideals of $A$.

An ideal $I$  of $A$ is called a  maximally differential ideal if
it is maximally $\mathfrak{D}$-differential for a set
$\mathfrak{D}$ of derivations or of higher derivations.

The  ring $A$ is called a differentially simple ring if the ideal
$(0)$ is maximally differential in $A$.

For a ring $A$ and an $A$-module $M$, let $\ell_A(M)$ denote the
length of $M$ over $A$.

We start by proving the following result:

\begin{lem}\label{l1}
Let $(A, {\frak m})$ be a Noetherian local ring and let
$\{x_1,x_2,\dots,x_n\}$ be a set of generators of ${\frak m}$. For
integers $k_1,\dots,k_n\geq 1${\rm ,}
$\ell_A(A/(x_1^{k_1},x_2^{k_2},\dots,x_n^{k_n}))\leq k_1\dots
k_n$.
\end{lem}

\begin{proof}
The proof is straightforward and perhaps it is available
somewhere. Unfortunately, we do not have a reference and therefore
we provide here a proof.

Note that the length of $A/(x_1^{k_1},x_2^{k_2},\dots,x_n^{k_n})$
is finite as its support is $\{{\frak m}\}$.

We prove the result by induction on $m={k_1}+{k_2}+\cdots+{k_n}$.
If $m=n$, then each $k_i$ equals 1 and the result is obvious. Now
suppose that $m> n$. Without loss of generality, we may assume
that $k_1\geq 2$. We then have the following exact sequence of
$A$-modules
\begin{align*}
%\hskip -2.5pc
\frac{A}{(x_1,x_2^{k_2},\dots,x_n^{k_n})}\!\stackrel{\phi}{\longrightarrow}\!
\frac{A}{(x_1^{k_1},x_2^{k_2},\dots,x_n^{k_n})}\!\stackrel{\pi}{\longrightarrow}\!
\frac{A}{(x_1^{k_1-1},x_2^{k_2},\dots,x_n^{k_n})}\!\longrightarrow\!
0,
\end{align*}
where $\phi(1+(x_1,x_2^{k_2},\dots,x_n^{k_n}))= x_1^{k_1-1}+
(x_1^{k_1},x_2^{k_2},\dots,x_n^{k_n})$ and $\pi$ is the natural
surjection. Hence
\begin{align*}
&\ell_A(A/(x_1^{k_1},x_2^{k_2},\dots,x_n^{k_n}))\\[.4pc]
&\quad\,\leq \ell_A(A/(x_1^{k_1-1},x_2^{k_2},\dots,x_n^{k_n}))
+ \ell_A(A/(x_1,x_2^{k_2},\dots,x_n^{k_n}))\\[.4pc]
&\quad\,\leq
(k_1-1){k_2}\dots{k_n}+{k_2}\dots{k_n}={k_1}{k_2}\dots{k_n}.
\end{align*}
\end{proof}

The following result is an extension of our result (Theorem~4 of
\cite{akm1}) in the sense that here we do not assume that the ring
contains a field and hence is valid even for local rings with
unequal characteristics.

\begin{theor}[\!]\label{pi}
Let $(A, {\frak m})$ be a Noetherian local ring and $\frak{D}$ be
a set of derivations of $A$. Let $I$ be a maximally differential
ideal in $A$. Let $n=\emdim(A)$ and $r=\emdim(A/I)$. If $A/I$ is
of characteristic $p>0${\rm ,} then there exists a minimal set
$\{x_1,x_2,\dots,x_n\}$ of generators of ${\mathfrak m}$ such that
$I = (x_1^p, \dots, x_r^p, x_{r+1}, \dots, x_n)$.
\end{theor}

\begin{proof}
For $d\in \frak{D}$, let $d'$ denote the derivation induced by $d$
on $A/I$ and let $\frak{D}'=\{d'|d\in \frak{D} \}$. Then $A/I$ is
differentially simple under $\frak{D}'$.

If $a\in\frak{m}$, then the ideal $(a^p)$ is
$\frak{D}$-differential. Therefore we have ${a}^p\in I$.

Choose a minimal set of generators $x_1,\dots,x_n$ of $\frak{m}$
such that $x_{r+1},\dots,x_n\in I$ and $x_1+I,\dots,x_r+I$ form a
minimal set of generators of $\frak{m}/I$. By Harper's theorem
(\cite{lh} or Corollary~3 of \cite{akm1}), $A/I$ contains a
coefficient field $k$ and the $k$-algebra map
\begin{equation*}
\phi\hbox{\rm :}\ k[y_1,\dots,y_r]\longrightarrow A/I,
\end{equation*}
defined by $\phi(y_i)=x_i+I$ (where $y_1,\dots,y_r$ are
indeterminates) for $i=1,\dots,r$; is onto with kernel
$(y_1^p,\dots,y_r^p)$. Therefore
$\ell_A(A/I)=\ell_{A/I}(A/I)=p^r$.

On the other hand,
$J=(x_1^p,\dots,x_r^p,x_{r+1},\dots,x_n)\subseteq I$ and therefore
$p^r=\ell_A(A/I)\leq\ell_A(A/J)\leq p^r$, where the last
inequality  follows by Lemma~\ref{l1}.

Hence $I=J=(x_1^p,\dots,x_r^p,x_{r+1},\dots,x_n)$.
\end{proof}

\begin{rem}
{\rm In view of the above result we see that if $(A,\frak{m})$ is
a regular local ring of dimension $n$, $I$ is an ideal of $A$
which is maximally differential under a set of derivations of $A$
and $A/I$ is of characteristic $p>0$. Then
$I=(x_1^p,\dots,x_r^p,x_{r+1},\dots,x_n),$ for some regular set
$x_1,\dots,x_n$ of parameters  of $A$ and therefore $\dim(A/I)=0$,
and hence in general, $A$ is not regular. Therefore the question
regarding the regularity of $A/I$ needs to be answered only when
\begin{enumerate}
\renewcommand\labelenumi{(\arabic{enumi})}
\leftskip .1pc
\item $A$ contains a field of characteristic zero and $I$ is
maximally differential under a set of derivations or

\item $I$ is maximally differential under a set of higher
derivations.\vspace{-.5pc}
\end{enumerate}}
\end{rem}
Also note that the first case is included in the second case.

We recall the following result from \cite{bs}:

\begin{theor}[(see Theorem~1.4 of \cite{bs})]\label{bst} Let $A$ be a
Noetherian local ring and let $I$ be a proper ideal of $A$.
Suppose that $I$ is maximally $\frak D$-differential for a set
$\frak D$ of higher derivations of $A$. Then $A$ is normally flat
along $I${\rm ,} that is{\rm ,} for all $n\geq 0${\rm ,}
$I^n/I^{n+1}$ is free as an $A/I$-module.
\end{theor}

We now prove our main result:

\begin{theor}[\!]
Let $A$ be a regular local ring and let $\mathfrak{D}$ be either a
set of derivations of $A$ or a set of higher derivations of $A$.
Let $I$ be the maximally $\mathfrak{D}$-differential ideal in $A$.
Then $I$ is generated by a regular sequence and hence $A/I$ is a
complete intersection.
\end{theor}

\begin{proof}
We need to show that $I$ is generated by an $A$-sequence. As $A$
is regular, by Auslander--Buchsbaum--Serre theorem (Theorem~2.2.7
of \cite{bh}), $A/I$ has a finite projective dimension over $A$.
Therefore, if we show that $I/I^2$ is free as an $A/I$-module,
then by Ferrand--Vasconcelos theorem (Theorem~2.2.8 of \cite{bh}),
it follows that $I$ is generated by an $A$-sequence.

First suppose that $\mathfrak D$ is a set of higher derivations of
$A$. By Theorem~\ref{bst}, $I/I^2$ is free as an $A/I$-module.
Therefore in this case we are through.

Now suppose that $\mathfrak{D}$ is a set of derivations of $A$.
For $d\in \mathfrak{D}$, let $d'$ denote the derivation induced by
$d$ on $A/I$ and let ${\mathfrak D}'=\{d'|d\in \mathfrak{D}\}$.
Then $A/I$ is differentially simple under the set $\mathfrak{D}'$.
Let $k=\{a\in A/I|d(a)=0$ for all $d \in \mathfrak{D}' \}$. By
differential simplicity of $A/I$, it follows that $k$ is a field.

If $A/I$ contains a field of characteristic zero then so does $A$
and hence $I$ is, in fact, maximally differential under the set
\begin{equation*}
\{(1,d,d^2/2!,\dots,d^n/n!,\dots)|d\in {\mathfrak D}\}
\end{equation*}
of higher derivations of $A$. Therefore, this case is included in
the previous case.

Now assume that $A/I$ is of prime characteristic $p>0$. The proof
that $I$ is generated by an $A$-sequence follows directly from
Theorem~\ref{pi}. By Theorem~\ref{pi}, there exists a minimal set
of generators $x_1,\dots,x_n$ of $\mathfrak m$ such that
$I=(x_1^p,\dots,x_r^p,x_{r+1},\dots,x_n),$ where
$n=\emdim(A)=\dim(A)$, $r=\emdim(A/I)$. As $A$ is regular,
$x_1,\dots,x_n$ forms an\break $A$-sequence, and therefore $I$ is
generated by an $A$-sequence.

That is, $A/I$ is a complete intersection.
\end{proof}

\begin{rem}
{\rm If $A$ is not regular then $A/I$ need not be a complete
intersection or even Cohen--Macaulay, as is clear from the example
constructed in \cite{bdl}.}
\end{rem}

\end{document}